\documentclass[12pt,twoside]{amsart}
\usepackage{amssymb}
\usepackage{amscd}

\title[Notes on Toric Varieties]{Notes on Toric Varieties from 
Mori theoretic viewpoint}
\author{Osamu Fujino}
\date{2001/12/8}
\address{Research Institute for Mathematical Sciences\\ 
Kyoto University, Kyoto 606-8502 Japan}
\email{fujino@kurims.kyoto-u.ac.jp}
\newcommand{\bQ}[0]{{\mathbb Q}}

\newcommand{\Supp}[0]{{\operatorname{Supp}}}
\newcommand{\xPic}[0]{{\operatorname{Pic}}}
\newcommand{\mult}[0]{{\operatorname{mult}}}

\newtheorem{thm}{Theorem}[section]
\newtheorem{lem}[thm]{Lemma}
\newtheorem{cor}[thm]{Corollary}
\newtheorem{prop}[thm]{Proposition}

\theoremstyle{definition}

\newtheorem{rem}[thm]{Remark}
\newtheorem*{ack}{Acknowledgements}       

\newtheorem*{co}{Comment}       

\newtheorem{say}[thm]{}
\newtheorem*{claim}{Claim}       

\theoremstyle{remark}

\newtheorem{step}{Step}
\begin{document}
\bibliographystyle{amsalpha+}
\maketitle

\abstract 
The main purpose of this notes is to supplement the paper 
\cite{reid}, which treated Minimal Model Program 
(also called Mori's Program) on toric 
varieties. 
We calculate lengths of negative extremal rays of 
toric varieties. As an application, we obtain 
a generalization of Fujita's conjecture for singular toric 
varieties. 
We also prove that 
every toric variety has a small projective toric 
$\bQ$-factorialization. 
\endabstract

\setcounter{section}{-1}
\section{Introduction}\label{se1}
The main purpose of this notes is to supplement the paper 
\cite{reid}, which treated Minimal Model Program 
(also called Mori's Program) on 
toric varieties. 
We calculate lengths of negative extremal rays of 
toric varieties. 
It is an easy exercise once we understand \cite {reid}. 
As a corollary, 
we obtain a strong version of Fujita's conjecture for 
singular toric varieties. 
Related topics are \cite {fujita}, \cite{kawamata}, \cite{l} and 
\cite[Section 4]{mu}. 
We will freely use 
the notation in \cite {fulton}, \cite{reid} and 
work over an algebraically closed field $k$ of arbitrary characteristic 
throughout this 
paper. 

\vspace{5mm}
The following is the main theorem of this paper. 

\begin{thm}[Cone Theorem]\label{main}
Let $X$ be an $n$-dimensional 
{\em{(}}not necessarily $\bQ$-factorial{\em{)}} 
projective toric variety over $k$. 
We can write the cone of curves as follows: 
$$NE(X)=\sum \mathbb R_{\geq 0}[C].$$ 
Let $D=\sum _j d_j D_j$ be a $\bQ$-divisor, 
where $D_j$ is an irreducible torus invariant divisor and 
$0\leq d_j\leq 1$ for every $j$.
Assume that $K_X+D$ is $\bQ$-Cartier. 
Then, for each extremal ray $\mathbb R_{\geq 0}[C]$, 
there exists an $(n-1)$-dimensional cone $\tau$ such that 
$[V(\tau)]\in \mathbb R_{>0}[C]$ and 
$$-(K_X+D)\cdot V(\tau)\leq n+1.$$ 
Moreover, we can choose 
$\tau$ such that $-(K_X+D)\cdot V(\tau)\leq n$ 
unless $X\simeq \mathbb P^{n}$ and $\sum_j d_j<1$. 
\end{thm}

Section \ref{sec1} deals with the proof of Theorem \ref{main}. 
We will treat an application of this theorem in Section \ref{2}. 
Professor Kajiwara informed me of \cite {mu} in Kinosaki 
after I finished 
the preliminary version of this paper. 
The following is Musta\c{t}\v{a}'s formulation of 
Fujita's conjecture for toric varieties. 
He proved it on the assumption that $X$ is 
non-singular as an application of his vanishing theorem (see 
\cite[Theorem 0.3] {mu}). 
Our proof doesn't use vanishing theorems. 

\begin{cor}[Strong version of Fujita's conjecture]\label{mufu}
Let $X$ be an $n$-dimensional {\em{(}}not necessarily 
$\bQ$-factorial{\em{)}} projective toric variety over $k$ and 
$D=\sum _j d_j D_j$ be a $\bQ$-divisor, 
where $D_j$ is an irreducible torus invariant divisor and 
$0\leq d_j\leq 1$ for every $j$. 
Assume that $K_X+D$ is $\bQ$-Cartier. 
Let $L$ be a line bundle on $X$. 
\begin{itemize}
\item[(1)] Suppose that 
$(L\cdot C)\geq n$ for every torus invariant integral 
curve $C\subset X$. Then $K_X+D+L$ is nef unless 
$X\simeq \mathbb P^{n}$, $\sum_j d_j<1$ and $L\simeq \mathcal O_{\mathbb P^n}
(n)$.

\item[(2)] Suppose that 
$(L\cdot C)\geq n+1$ for every torus invariant integral 
curve $C\subset X$. Then $K_X+D+L$ is ample unless 
$X\simeq \mathbb P^{n}$, $D=0$ and $L\simeq \mathcal O_{\mathbb P^n}
(n+1)$. 
\end{itemize}
\end{cor}
Of course, we can recover \cite[Theorem 0.3]{mu} easily if we 
assume that $X$ is non-singular. See also Remark \ref{remrem}.  

In section \ref{MMP}, we collect some results obtained by 
Minimal Model Program on toric varieties. 
We need Lemma \ref{ro} for the proof of Theorem \ref{main}. 
We prove that 
every toric variety has a small projective toric 
$\bQ$-factorialization. For related topics, see \cite 
[Section 3]{op}. 

\begin{co}
After I circulated \cite{f}, I found \cite{l}. 
In \cite{l}, Laterveer proved Fujita's conjecture for 
$\bQ$-Gorenstein projective toric varieties 
by the similar argument to mine. 
\end{co}

\begin{ack}Some parts of this paper were obtained 
in 1999, when I was a Research Fellow of the Japan Society for 
the Promotion of Science. I would like 
to express my gratitude to Professors 
Masanori Ishida, Shigefumi Mori, Tadao Oda, 
Takeshi Kajiwara, and Hiromichi Takagi, who gave me various advice 
and useful comments. I like to thank Doctor Hiroshi Sato, 
who gave me various advice and answered my questions. 
I also like to thank Doctor Takeshi Abe, who led me to this problem. 
\end{ack}

\section{Proof of the theorem}\label{sec1}
\begin{say} 
First, let us recall {\em{weighted projective spaces}}. 
We adopt toric geometric descriptions. 
This helps the readers to understand Theorem \ref{main}. 
However, it is not necessary for the proof of Theorem \ref{main}. 
\end{say}

\begin{say}[c.f. {\cite[p.35]{fulton}}]\label{1} 
Let $\mathbb P(d_1,\cdots, d_{n+1})$ 
be a {\em{weighted projective space}}. 
To construct this as  a toric variety, 
start with the fan whose cones 
generated by proper subsets 
of $\{v_1,\cdots, v_{n+1}\}$, where any $n$ of these vectors 
are linearly independent, 
and their sum is zero. 
The lattice $N$ is taken to be generated by the vectors 
$e_i=(1/{d_i})\cdot v_i$ for $1\leq i \leq n+1$. 
The resulting toric variety is in fact 
$\mathbb P=\mathbb P(d_1,\cdots, d_{n+1})$. 
We note that $\xPic \,\mathbb P\simeq \mathbb Z$. 
Let $f_i$ be a unique primitive lattice point in the 
cone $\langle e_i\rangle$ with 
$e_i=u_i f_i$ for $u_i\in \mathbb Z_{>0}$. 
We put $d=\gcd (u_1 d_1,\cdots, u_{n+1}d_{n+1})$ 
and define $c_i=(1/d)u_id_i$ for every $i$. 
Then we obtain that $\mathbb P(d_1,\cdots, d_{n+1})\simeq 
\mathbb P(c_1,\cdots, c_{n+1})$ and 
$\sum c_i f_i=0$. 
By changing the order, we can assume that 
$c_1\leq c_2 \leq \cdots \leq c_{n+1}$. 
We note that $-K_{\mathbb P}=\sum V(f_i)$. 
Let $\tau$ be the $(n-1)$-dimensional cone 
$\langle f_1,\cdots, f_{n-1}\rangle$. 
Then we have $$-K_{\mathbb P} \cdot V(\tau)=
\sum _{i=1}^{n+1} V({f_i})\cdot V(\tau) =
\frac {c}{c_{n}c_{n+1}}(\sum_{i=1}^{n+1} c_i)\leq n+1, $$
where $c=\gcd (c_n, c_{n+1})$. 
We note that 
$$V(f_i)\cdot V(\tau)=\frac{c\ c_i}{c_n c_{n+1}}\ .$$
For calculations of 
intersection numbers, 
we recommend the readers to 
see \cite[p.100]{fulton} and \cite[(2.7)]{reid}. 
If the equality holds in the above equation, 
then $c_i=1$ for every $i$. 
Thus, we obtain $\mathbb P\simeq \mathbb P^n$. 

\begin{rem}\label{pripri}
Suppose that $\gcd (d_1,\cdots,d_{n+1})=1$. Then, 
$e_i$ is primitive in $\langle e_i\rangle\cap N$ 
if and only if $\gcd (d_1,\cdots, d_{i-1},d_{i+1},
\cdots, d_{n+1})=1$. 
\end{rem}
\end{say}

\begin{say}\label{nijigen}
If $n=2$, then $c=1$ since $f_1$ is primitive and 
$\sum c_i f_i=0$. 
Therefore, we have 
$$-K_{\mathbb P} \cdot V(\tau)=
\frac {1}{c_{2}c_{3}}(\sum_{i=1}^{3} c_i)\leq 
\frac{1}{2}+\frac {1}{2}+1\leq 2=n $$
when $\mathbb P\not\simeq \mathbb P^2$. 
So, we have that $-K_{\mathbb P} \cdot V(\tau)\leq n$ if 
$n=2$ and $\mathbb P\not\simeq \mathbb P^2$. 
If $-K_{\mathbb P} \cdot V(\tau)=2$, then 
$\mathbb P\simeq \mathbb P(1,1,2)$. 
\end{say}

\begin{say}[c.f. Proposition \ref{key} below]\label{dame}
When $n\geq 3$, the above inequality in \ref {nijigen} is 
not true. Assume that $n\geq 3$. 
Let $\mathbb P$ be an 
$n$-dimensional weighted projective space $\mathbb P(l-1,l-1,l,\cdots,l)$, 
where $l\geq 2$. 
Then we obtain $$-K_{\mathbb P} \cdot V(\tau)=n+1-\frac{2}{l}.$$ 
So, we have $-K_{\mathbb P} \cdot V(\tau)>n$ when $l\geq 3$. 
If we make $l$ large, then $-K_{\mathbb P} \cdot V(\tau)$ becomes 
close to $n+1$. 
\end{say}

\begin{say}
Let $\mathbb P=
\mathbb P(1,\cdots,1,l-1,l)$ be an $n$-dimensional 
weighted projective 
space with $l\geq 2$ and 
$n\geq 2$. Then we have 
$$
-K_{\mathbb P} \cdot V(\tau)=\frac{n+2l-2}{l(l-1)}. 
$$ 
Thus, if we make $l$ large, then 
$-K_{\mathbb P} \cdot V(\tau)$ becomes close to 
zero.  
\end{say}

\begin{say}
Next, we treat $\bQ$-factorial 
toric Fano varieties with Picard number 
one. This type of varieties plays an important 
role for the analysis of extremal contractions. 
Here, we adopt the following description \ref{tuitui} for 
the definition of $\bQ$-factorial 
toric Fano varieties with Picard number 
one. By this, it is easy to see that 
every extremal contraction contains them 
in the fibers (see Proof of the theorem below). 
Of course, weighted projective spaces are in this class. 

\begin{rem}
In \cite[(0.1)]{reid}, it is stated that any fiber of 
an extremal contraction is a weighted projective space. 
However, it is not true since there exists a 
$\bQ$-factorial toric Fano variety with Picard number one 
that is not a weighted projective space.  
\end{rem}
\end{say}

\begin{say}[$\bQ$-factorial 
toric Fano varieties with Picard number one]\label{tuitui}
Now we fix $N\simeq \mathbb Z ^n$. Let $\{v_1,\cdots,v_{n+1}\}$ 
be a set of primitive vectors such that $N_{\mathbb R}=\sum _i 
\mathbb R_{\geq 0}v_i$. 
We define $n$-dimensional cones 
$$
\sigma_i:=\langle e_1,\cdots,e_{i-1},e_{i+1},\cdots,e_{n+1}\rangle 
$$ 
for $1\leq i\leq n+1$. 
Let $\Delta$ be the complete fan generated by $n$-dimensional 
cones $\sigma_i$ for every $i$. Then 
we obtain a complete toric variety $X=X(\Delta)$ with 
Picard number $\rho (X)=1$. 
We call it a {\em{$\bQ$-factorial 
toric Fano variety with Picard number one}}. 
We define $(n-1)$-dimensional cones $\mu_{i,j}=\sigma _i\cap \sigma _j$ 
for $i\ne j$. 
We can write $\sum _i a_i v_i=0$, where $a_i\in \mathbb Z_{>0}$, 
$\gcd(a_1,\cdots,a_{n+1})=1$, and $a_1\leq a_2\leq\cdots\leq a_{n+1}$ 
by changing the order. 
Then we obtain 
$$
0< V({v_{n+1}})\cdot V(\mu_{n,n+1})=\frac{\mult {(\mu_{n,n+1})}}
{\mult {(\sigma_{n})}}\leq 1, 
$$
$$
V({v_{i}})\cdot V(\mu_{n,n+1})=\frac{a_i}{a_{n+1}}\cdot
\frac{\mult {(\mu_{n,n+1})}}
{\mult {(\sigma_{n})}}, 
$$
and 
\begin{eqnarray*}
-K_{X} \cdot V(\mu_{n,n+1})&=&
\sum _{i=1}^{n+1} V({v_i})\cdot V(\mu_{n,n+1})\\
& =&
\frac {1}{a_{n+1}}
{(\sum_{i=1}^{n+1} a_i)}
\frac{\mult {(\mu_{n,n+1})}}
{\mult {(\sigma_{n})}}\leq n+1. 
\end{eqnarray*} 
For ``$\mult$'' in the above equations, 
see \cite[p.48 and p.100]{fulton}. 
If $-K_{X} \cdot V(\mu_{n,n+1})=n+1$, then $a_i=1$ for every 
$i$ and $\mult (\mu_{n,n+1})=\mult (\sigma_{n})$. 

\begin{prop}\label{key}
If $X\not\simeq \mathbb P^n$, then there exists some pair 
$(l,m)$ such that $-K_{X} \cdot V(\mu_{l,m})\leq n$. 
\end{prop}
\proof
Assume the contrary. 
Then we obtain 
$$-K_{X} \cdot V(\mu_{k,n+1})=\frac {1}{a_{n+1}}
{(\sum_{i=1}^{n+1} a_i)}
\frac{\mult {(\mu_{k,n+1})}}
{\mult {(\sigma_{k})}}>n
$$ 
for $1\leq k\leq n$. 
Thus 
$$
(n+1)a_{n+1}\geq \sum_{i=1}^{n+1} a_i 
>\frac{\mult {(\sigma_{k})}}
{\mult {(\mu_{k,n+1})}}na_{n+1}  
$$
for every $k$. 
Since 
$$
\frac{\mult {(\sigma_{k})}}
{\mult {(\mu_{k,n+1})}}\in \mathbb Z_{>0}, 
$$
we have that $\mult {(\sigma_{k})}=\mult (\mu_{k,n+1})$ 
for every $k$. 
This implies that $a_{k}$ divides $a_{n+1}$ for all $k$. 
\begin{claim}
$a_1=\cdots=a_{n+1}=1$. 
\end{claim}
\proof[Proof of Claim]
If $a_1=a_{n+1}$, then we obtain the required results. 
So, we assume that $a_1\ne a_{n+1}$. 
On this assumption, we have that $a_2\ne a_{n+1}$ since 
$v_1$ is primitive and $\sum _i a_i v_i=0$. 
In this case, we have 
$$-K_{X} \cdot V(\mu_{k,n+1})=\frac {1}{a_{n+1}}
{(\sum_{i=1}^{n+1} a_i)}\leq n. 
$$ 
We note that 
$$
\frac{a_i}{a_{n+1}}\leq \frac{1}{2} 
$$ 
for $i=1,2$. This is a contradiction. So we obtain 
that $a_1=\cdots=a_{n+1}=1$. 
\endproof
In this case, $-K_X\cdot V(\mu_{i,j})>n$ implies 
$-K_X\cdot V(\mu_{i,j})=n+1$ for every pair $(i,j)$. 
Then $\mult (\mu_{i,j})=\mult (\sigma_{i})$ for 
$i\ne j$. 
So, we have that $\mult (\sigma_{i})=1$ for every $i$. 
Therefore, we obtain $X\simeq \mathbb P^n$. 
This is a contradiction. 
\endproof

\begin{rem}
The usual definition of Fano varieties is 
the following: $X$ is Fano if $-K_X$ is an ample $\bQ$-Cartier divisor. 
It is easy to check that 
the notion of 
$\bQ$-factorial toric Fano varieties with Picard number one 
by the usual definition coincides with 
ours. 
\end{rem}
\end{say}

\begin{say}
From now on, we freely use the notation in \cite {reid}, 
especially, \cite [(2.2)]{reid} (see also \cite[(1.10)]{reid}).  

\proof[Proof of the theorem] 
\begin{step}\label{st} 
We assume that $X$ is $\bQ$-factorial. 
Let $R=\mathbb R_{\geq 0}[C]$ be an extremal ray. 
There exists an elementary contraction 
$\varphi _R:X\to Y$, 
which is corresponding to the extremal ray $R$. 
The $\bQ$-factorial 
toric Fano variety $P\subset X$ with Picard number $\rho (P)=1$, 
which is corresponding to the cone 
$\sigma=\langle e_1,\cdots,e_{\beta}
\rangle$, 
that is, $P=V(\sigma)\subset X$, is 
a fiber of $\varphi _R|_{A}:A\to B$ (c.f. \cite[(2.5)]{reid}).  
We note that 
$$K_P=-\sum _{i=\beta +1}^{n+1}V(\widetilde \rho _i),$$ 
where $\widetilde \rho _i=\langle e_1,
\cdots,e_{\beta}, e_i\rangle$ 
for $\beta +1\leq i \leq n+1$. 
On the other hand, $V(\widetilde \rho _i)=b_iV(e_i)\cdot V(\sigma)$ 
for some $b_i\in \mathbb Z_{>0}$ since the cones are simplicial 
(see \cite[p.100]{fulton}). 
Let $\widetilde \tau$ be an $(n-1)$-dimensional cone containing 
$\sigma$. 
We have that 
\begin{eqnarray*}
K_P\cdot V(\widetilde \tau)&
=&-\sum _{i=\beta+1}^{n+1}V(\widetilde \rho _i)
\cdot V(\widetilde \tau)\\&=&
-V(\widetilde \tau)\cdot
(\sum _{i=\beta+1}^{n+1}b_iV(e_i)\cdot V(\sigma))
\\&=&V(\widetilde \tau)\cdot
(K_X+\sum_{\text{every ray}}V(e_i)-\sum _{i=\beta+1}^{n+1}
b_iV(e_i))
\\&=&V(\widetilde \tau)\cdot(K_X+\sum _{i=\beta+1}^{n+1}
(1-b_i)V(e_i)+\sum_{\text{others}}V(e_i))
\\&\leq& (K_X+D)\cdot V(\widetilde \tau). 
\end{eqnarray*}
We note that $$K_X+\sum_{\text{every ray}}V(e_i)\sim 0$$ 
and $D$ can be written as $\sum_j d_j V(e_j)$ 
with $0\leq d_j\leq 1$ by the assumption, and 
that $ V(\widetilde \tau)\cdot
V(e_i)>0$ if and only if $\beta+1\leq i\leq n+1$ 
by \cite[(2.2)]{reid} (see also \cite [(2.4), (2.7), (2.10)]{reid}). 
We choose $\widetilde \tau$ as in the above argument \ref{tuitui}, 
that is, $-K_P\cdot V(\widetilde \tau) 
\leq n-\beta+1$, where $\dim P=n-\beta$. 
Then, by the above argument and the choice of $\widetilde \tau$, 
$$-(K_X+D)\cdot V(\widetilde \tau)
\leq -K_P\cdot V(\widetilde \tau)\leq 
n-\beta+1. 
$$ 
Therefore, if the minimal length of a 
$(K_X+D)$-negative extremal ray is 
greater than $n$, then $\beta=\alpha =0$. 
Thus we have $X\simeq \mathbb P^n$ and $\sum_j d_j <1$ 
by Proposition \ref{key}. 
Therefore, we obtain the required result when 
$X$ is $\bQ$-factorial. 
\end{step}
\begin{step}[{c.f. \cite[(2.4) Lemma]{l}}]
We assume that $X$ is not $\bQ$-factorial. 
Let $f:(\widetilde X, \widetilde D)\to (X,D)$ be a 
projective modification constructed in 
Lemma \ref{ro} below. 
We note that $X\not\simeq \mathbb P^n$. 
Let $R=\mathbb R_{\geq 0}[C]$ be a 
$(K_X+D)$-negative extremal ray. 
We take $V(\tau)\in \mathbb R_{> 0}[C]$ such that 
$-(K_X+D)\cdot V(\tau)$ is minimal. 
We take $V(\widetilde \tau)$ on $\widetilde X$ 
such that $f_*V(\widetilde \tau)=V(\tau)$. 
We can write $V(\widetilde\tau)=\sum a_i
V(\widetilde\tau_i)$ in $NE(\widetilde X)$ for 
$a_i\in \mathbb R_{>0}$ such that $V(\widetilde \tau_i)$ 
is extremal and $-(K_{\widetilde X}+\widetilde D)\cdot 
V(\widetilde \tau _i)\leq n$ for every $i$ 
by Theorem \ref {main} 
since $\widetilde X$ is not a projective space. 
Since $\sum_i a_i f_*V(\widetilde \tau _i)=V(\tau)\in R$, 
we have that $f_*V(\widetilde \tau _i)\in R$ for every $i$. 
So, there exists some $i$ such that 
$0\ne f_*V(\widetilde \tau _i)=bV(\tau)$ in $R$ for $b\geq 1$ 
since $-(K_X+D)\cdot V(\tau)$ is minimal. 
Therefore, 
$$-(K_X+D)\cdot V(\tau)=-\frac{1}{b}(K_{\widetilde X}+
\widetilde D)\cdot V(\widetilde\tau_i)\leq n. 
$$ 
Thus we finished the proof. 
\endproof
\end{step}

\begin{rem}\label{remrem}
In Step \ref{st} in 
the proof of the theorem, we assume that $X$ is non-singular. 
Then we obtain that $b_i=1$ and 
$V(\widetilde \tau)\cdot V(e_i)\in \mathbb 
Z$. We note that $V(\widetilde \tau)\cdot V(e_i)>0$ if and only if 
$\beta+1\leq i\leq n+1$. 
It can be checked easily 
that $P$ is an $(n-\beta)$-dimensional projective space $\mathbb 
P^{n-\beta}$ and 
$K_P\cdot V(\widetilde \tau)=-(n-\beta+1)$. 
Thus, Proposition 4.3, 
Lemma 4.4, and Propositions 4.5, 4.6 in \cite {mu} 
can be checked easily by the above computation 
(see also \cite[(2.10)(i)]{reid}). 
Therefore, we can recover \cite[Section 4]{mu} without 
using vanishing theorems.  
\end{rem}
\end{say}

\section{Applications to Fujita's conjecture on toric varieties}\label{2}

In this section, we treat some applications of 
Theorem \ref{main}. Corollary \ref{mufu} follows from 
Theorem \ref{main} directly. 

\proof[Proof of Corollary \ref{mufu}]
It is obvious by 
Theorem \ref{main}, \cite[\S2.3 Theorem 2.18] {oda} 
and \cite[Theorem 3.2]{mu}. 
\endproof

\begin{cor}
In Corollary \ref{mufu} (1), we further assume that 
$K_X+D$ is Cartier. 
Then $K_X+D+L$ is generated by global sections 
unless $X\simeq \mathbb P^n$, $D=0$, and 
$L\simeq \mathcal O_{\mathbb P^n}(n)$. 
\end{cor}
\proof
It is obvious by Corollay \ref{mufu} (1). 
We note that every nef line bundle is generated by its 
global sections on a complete toric variety. 
It is well-known (see, for example, \cite[Theorem 3.1]{mu}). 
\endproof

By combining Corollary \ref{mufu} with 
Demazur's theorem (\cite[\S2.3 Corollary 2.15]{oda}), 
we obtain the following result. This is 
the original version of Fujita's conjecture on 
toric varieties.   

\begin{cor}[Fujita's conjecture for toric varieties] 
Let $X$ be a non-singular projective toric variety over $k$ and 
$L$ an ample line bundle on $X$. 
Then $K_X+(n+1)L$ is generated by global 
sections and $K_X+(n+2)L$ is very ample, 
where $n=\dim X$. 
Moreover, if $(X,L)\not\simeq (\mathbb P^n, \mathcal O_{\mathbb P^n}(1))$, 
then $K_X+nL$ is generated by global 
sections and $K_X+(n+1)L$ is very ample. 
\end{cor}

\begin{rem}
For very ampleness on singular toric varieties, 
see \cite[3. Very ampleness]{l}. 
\end{rem}

\section{Remarks on Minimal Model Program for toric varieties}\label{MMP}

In this section, we use the notation in \cite {km} and \cite {fa}. 
The following lemma is well-known to specialists. 
It may help the readers to understand this section.  

\begin{lem}\label{sing}
Let $X$ be a complete toric variety over $k$ 
and $D$ the complement of the big torus in $X$ 
as a reduced divisor. 
Then the pair $(X,D)$ is log-canonical. 
Furthermore, if $K_X$ is $\bQ$-Cartier, then 
the pair $(X,0)$ is log-terminal. 
\end{lem}
\proof
Let $g:Y\to X$ be a toric resolution of singularities. 
Then  we have 
$$
K_Y+E=g^*(K_X+D), 
$$ 
where $E$ is the complement of the big torus 
in $Y$ as a reduced divisor. 
Thus, the pair $(X,D)$ is log-canonical by 
the definition (see \cite[Definition 2.34]{km}). 
If $K_X$ is $\bQ$-Cartier, then $D$ is 
$\bQ$-Cartier since $K_X+D\sim 0$. 
Note that $\Supp g^*D=\Supp E$ and $g^*D$ is 
an effective $\bQ$-divisor. 
Therefore, the pair $(X,0)$ is log-terminal (see 
\cite [Definition 2.34]{km}). 
\endproof

The following is a variant of \cite [17.10 Theorem]{fa} 
for toric varieties. We recommend the readers who are not familiar with 
Minimal Model Program to see \cite[\S 3.7]{km}. 

\begin{thm}\label{taisetsu}
Let $X$ be a complete toric variety over $k$ and 
$g:Y\to X$ a projective birational toric morphism from 
a $\bQ$-factorial toric variety $Y$. 
Let $\mathcal E$ be a subset of the exceptional divisors. 
Then there is a factorization 
$$
g:Y\dashrightarrow \widetilde X \to X
$$
with the following properties: 
\begin{itemize}
\item[(1)] $h:Y\dashrightarrow \widetilde X$ 
is a local isomorphism at every generic point 
of the divisor that is not in $\mathcal E$; 
\item[(2)] $h$ contracts every exceptional 
divisor in $\mathcal E$; 
\item[(3)] $\widetilde X$ is projective over $X$ and $\bQ$-factorial. 
Of course, the pair $(\widetilde X,0)$ is log-terminal by 
Lemma \ref{sing}.  
\end{itemize}
In particular, if $\mathcal E$ is the set of all the 
$g$-exceptional divisors, then $f:\widetilde X\to 
X$ is small, that is, an isomorphism in codimension one.   
We call this {\em{a small projective toric $\bQ$-factorialization}}. 
\end{thm}

\proof 
Let $g:Y\to X$ be as above and 
$E=\sum E_i$ the complement of the big torus 
in $Y$. 
We note that 
$$
K_Y+E=g^*(K_X+D)\sim 0. 
$$ 
Apply $(K_Y+\sum_{E_i\not \in \mathcal E} 
E_i+\sum _{E_j \in \mathcal E}
2E_j)
$-log minimal model program over $X$. 
We note that divisorial contractions and log-flips 
always exist by \cite[(0.1)]{reid} (see also \cite[\S 5-2]{kmm}). 
Here, a log-flip means an {\em{elementary transformation}} with respect 
to a $(K_Y+\sum_{E_i\not \in \mathcal E} E_i+\sum _{E_j \in \mathcal E}
2E_j)$-negative extremal ray in the terminology of \cite{reid}. 
Since the relative Picard number 
$\rho (Y/X)$ is finite, divisorial contractions can occur finite times. 
So it is enough to check the termination of log-flips. 
Assume that there exists an infinite sequence of log-flips: 
$$
Y_0\dashrightarrow Y_1\dashrightarrow \cdots \dashrightarrow Y_m
\dashrightarrow \cdots. 
$$
Let $\Delta$ be the fan corresponding to $Y_0$. 
Since the log-flips don't change one-dimensional 
cones of $\Delta$, there are numbers $k<l$ 
such that $Y_k\simeq Y_l$ over 
$X$. 
This is a contradiction because 
there is a valuation $v$ such that the discrepancies satisfy 
$$
a(v,Y_k,\sum_{E_i\not \in \mathcal E} 
E_i+\sum _{E_j \in \mathcal E}
2E_j)<a(v,Y_l,\sum_{E_i\not \in \mathcal E} 
E_i+\sum _{E_j \in \mathcal E}
2E_j)
$$ (see \cite[Lemma 3.38]{km}), where 
$\sum_{E_i\not \in \mathcal E} E_i+\sum _{E_j \in \mathcal E}2E_j$ 
means the proper transform of it on $Y_k$ or $Y_l$. 
Therefore, we obtain $f:\widetilde X\to X$ with the above mentioned 
properties by \cite[Lemma 3.39]{km}. 
\endproof

\begin{rem}\label{new}
Since we can take a projective toric desingularization 
as $g:Y\to X$ in Theorem \ref{taisetsu}, 
there exists at least one small projective toric 
$\bQ$-factorialization for $X$. 
\end{rem}

\begin{rem}[{c.f. \cite[Theorem 6.38]{km}}] 
Let $X$ be a complete toric variety and 
$f_i:X_i\to X$ be small projective toric 
$\bQ$-factorializations for $i=1,2$. 
Then $X_1$ and $X_2$ can be obtained from each 
other by a finite succession of {\em{elementary 
transformations}}\begin{footnote}
{This {\em{elementary transformation}} 
was called {\em{flop}} in \cite{op} (see 
\cite[p.397 Remark]{op}). However, it might be better 
to call it {\em{log-canonical flop}} from 
the log Minimal Model Theoretic viewpoint 
(c.f. Lemma \ref{sing}). See also 
\cite[6.8 Definition]{fa}.}.\end{footnote} 
It can be checked by the log Minimal Model Program over $X$ 
as in the proof of 
Theorem \ref{taisetsu} and \cite [Lemma 6.39]{km}. 
Details are left to the readers.  
\end{rem}

By Theorem \ref{taisetsu}, we obtain the next lemma, which 
was already used in the proof of Corollary \ref{mufu}.  

\begin{lem}\label{ro}
Let $X$ be a projective toric variety over $k$ 
and $D=\sum _j d_j D_j$ be a $\bQ$-divisor, 
where $D_j$ is an irreducible torus invariant divisor and 
$0\leq d_j\leq 1$ for every $j$. 
Assume that $K_X+D$ is $\bQ$-Cartier. 
Then there exists a projective birational toric morphism 
$f:\widetilde X\to X$ such that 
$\widetilde X$ has only $\bQ$-factorial singularities 
and $K_{\widetilde X}+\widetilde D=f^*(K_X+D)$, 
where $\widetilde D=\sum _i \widetilde d _i\widetilde D_i$ is a 
$\bQ$-divisor such that $\widetilde 
D_i$ is an irreducible torus invariant divisor and 
$0\leq \widetilde d_i\leq 1$ for every $i$. 
\end{lem}

By Sumihiro's equivariant embedding theorem, 
we can remove the assumption that 
$X$ is complete. 

\begin{cor}[Small projective toric $\bQ$-factorialization] 
Let $X$ be a toric variety over $k$. 
Then there exists a small projective toric morphism 
$f:\widetilde X\to X$ such that 
$\widetilde X$ is $\bQ$-factorial. 
\end{cor}
\proof
We can compactify $X$ by Sumihiro's theorem \cite[\S 1.4]{oda}. 
So, this corollary follows from Theorem \ref{taisetsu} 
and Remark \ref{new} easily. 
\endproof

The existence of a small projective toric $\bQ$-factorialization 
means the following corollary. 

\begin{cor}\label{cone}
Let $\Delta$ be a fan. 
Then there exists a projective simplicial subdivision $\widetilde
\Delta$ of $\Delta$, 
that is, the morphism $X(\widetilde \Delta)\to 
X(\Delta)$ is projective and $X(\widetilde \Delta)$ is 
$\bQ$-factorial, such that the set of one-dimensional 
cones of $\widetilde \Delta$ coincides with that of 
$\Delta$. 
\end{cor}

\begin{rem} 
The above corollary seems to follow from the 
theory of Gelfand-Kapranov-Zelevinskij decompositions. 
For details about GKZ-decompositions, 
see \cite[Section 3]{op}, especially, \cite[Corollary 3.8]{op}. 
We note that \cite{op} generalized and reformulated 
results on \cite{reid}. 
\end{rem}

\ifx\undefined\bysame
\newcommand{\bysame}{\leavevmode\hbox to3em{\hrulefill}\,}
\fi


\begin{thebibliography}{KMM}

\bibitem[F]{f} 
O.~Fujino, 
Notes on toric varieties from Mori theoretic 
viewpoint, RIMS-{\textbf {1340}} (2001). 

\bibitem[Ft]{fujita}
T.~Fujita, 
On polarized manifolds whose adjoint bundles are not semipositive, 
Algebraic geometry,
Sendai, 1985, 167--178, Adv. Stud. Pure Math., 
{\textbf{10}}, North-Holland, Amsterdam, 1987. 

\bibitem[Fl]{fulton}
W.~Fulton, 
{\em{Introduction to toric varieties}}, Annals 
of Mathematics Studies, {\textbf{131}}. The William H. Roever
Lectures in Geometry. Princeton University Press, Princeton, NJ, 1993. 

\bibitem[Ka]{kawamata}
Y.~Kawamata, 
On the length of an extremal rational curve, 
Invent. Math. {\textbf{105}} (1991), no. 3, 609--611. 

\bibitem[KMM]{kmm}
Y.~Kawamata, K.~Matsuda, and K.~Matsuki, Introduction to the Minimal 
Model Problem, in {\em Algebraic Geometry, Sendai 1985,} Advanced Studies 
in Pure Math. {\textbf {10}}, (1987) Kinokuniya and North-Holland, 283--360. 

\bibitem[KM]{km}
J.~Koll\'ar and S.~Mori, {\em{Birational geometry of 
algebraic varieties,}} Cambridge University Press (1998). 

\bibitem[L]{l} 
R.~Laterveer, 
Linear systems on toric varieties, 
T\^{o}hoku 
Math. J. {\textbf {48}}, (1996), 451--458.  

\bibitem[Mu]{mu}
M.~Musta\c{t}\v{a}, 
Vanishing Theorems on Toric Varieties, to appear in T\^{o}hoku 
Math. J. (math.AG/0001142). 

\bibitem[Od]{oda}
T.~Oda, 
{\em{Convex bodies and algebraic geometry}}. 
An introduction to the theory of toric varieties, 
Translated from the 
Japanese. Ergebnisse der Mathematik und ihrer Grenzgebiete (3) 
[Results in Mathematics and Related
Areas (3)], {\textbf{15}}. Springer-Verlag, Berlin, 1988.

\bibitem[OP]{op} 
T.~Oda and H.~S.~Park, Linear Gale transforms 
and Gelfand-Kapranov-Zelevinskij decompositions, 
T\^ohoku Math.~J. {\textbf{43}} (1991), 357--399. 

\bibitem[Re]{reid}
M.~Reid, Decomposition of toric morphisms, 
Arithmetic and geometry, Vol. II, 395--418, Progr. Math., {\textbf{36}},
Birkh\"auser Boston, Boston, MA, 1983. 

\bibitem[Ut]{fa}
J.~Koll\'ar, et al, {\em Flips and Abundance for 
Algebraic Threefolds,} Ast\'erisque {\textbf {211}}, Soc. Math. de. France, 
1992. 
\end{thebibliography}
\end{document}